\documentclass{article}
\usepackage{amssymb,amsmath,amsthm}

\newcommand{\T}{\mathbf{\sf L}\colon}
\newcommand{\F}{\mathbf{\sf R}\colon}

\newtheorem{definition}{Definition}

\newtheorem{theorem}{Theorem}
\newtheorem{proposition}[theorem]{Proposition}
\newtheorem{corollary}[theorem]{Corollary}
\newtheorem{lemma}[theorem]{Lemma}

\makeatletter
\newdimen\proofrulebreadth \proofrulebreadth=.05em
\newdimen\proofdotseparation \proofdotseparation=1.25ex
\newdimen\proofrulebaseline \proofrulebaseline=2ex
\newcount\proofdotnumber \proofdotnumber=3
\let\then\relax
\def\hfi{\hskip0pt plus.0001fil}
\mathchardef\squigto="3A3B
\newif\ifinsideprooftree\insideprooftreefalse
\newif\ifonleftofproofrule\onleftofproofrulefalse
\newif\ifproofdots\proofdotsfalse
\newif\ifdoubleproof\doubleprooffalse
\let\wereinproofbit\relax
\newdimen\shortenproofleft
\newdimen\shortenproofright
\newdimen\proofbelowshift
\newbox\proofabove
\newbox\proofbelow
\newbox\proofrulename
\def\shiftproofbelow{\let\next\relax\afterassignment\setshiftproofbelow\dimen0 }
\def\shiftproofbelowneg{\def\next{\multiply\dimen0 by-1 }%
\afterassignment\setshiftproofbelow\dimen0 }
\def\setshiftproofbelow{\next\proofbelowshift=\dimen0 }
\def\setproofrulebreadth{\proofrulebreadth}

\def\prooftree{
\ifnum  \lastpenalty=1
\then   \unpenalty
\else   \onleftofproofrulefalse
\fi
\ifonleftofproofrule
\else   \ifinsideprooftree
        \then   \hskip.5em plus1fil
        \fi
\fi
\bgroup
\setbox\proofbelow=\hbox{}\setbox\proofrulename=\hbox{}%
\let\justifies\proofover\let\leadsto\proofoverdots\let\Justifies\proofoverdbl
\let\using\proofusing\let\[\prooftree
\ifinsideprooftree\let\]\endprooftree\fi
\proofdotsfalse\doubleprooffalse
\let\thickness\setproofrulebreadth
\let\shiftright\shiftproofbelow \let\shift\shiftproofbelow
\let\shiftleft\shiftproofbelowneg
\let\ifwasinsideprooftree\ifinsideprooftree
\insideprooftreetrue
\setbox\proofabove=\hbox\bgroup$\displaystyle 
\let\wereinproofbit\prooftree
\shortenproofleft=0pt \shortenproofright=0pt \proofbelowshift=0pt
\onleftofproofruletrue\penalty1
}

\def\eproofbit{
\ifx    \wereinproofbit\prooftree
\then   \ifcase \lastpenalty
        \then   \shortenproofright=0pt  
        \or     \unpenalty\hfil         
        \or     \unpenalty\unskip       
        \else   \shortenproofright=0pt  
        \fi
\fi
\global\dimen0=\shortenproofleft
\global\dimen1=\shortenproofright
\global\dimen2=\proofrulebreadth
\global\dimen3=\proofbelowshift
\global\dimen4=\proofdotseparation
\global\count255=\proofdotnumber
$\egroup  
\shortenproofleft=\dimen0
\shortenproofright=\dimen1
\proofrulebreadth=\dimen2
\proofbelowshift=\dimen3
\proofdotseparation=\dimen4
\proofdotnumber=\count255
}

\def\proofover{
\eproofbit 
\setbox\proofbelow=\hbox\bgroup 
\let\wereinproofbit\proofover
$\displaystyle
}%
\def\proofoverdbl{
\eproofbit 
\doubleprooftrue
\setbox\proofbelow=\hbox\bgroup 
\let\wereinproofbit\proofoverdbl
$\displaystyle
}%
\def\proofoverdots{
\eproofbit 
\proofdotstrue
\setbox\proofbelow=\hbox\bgroup 
\let\wereinproofbit\proofoverdots
$\displaystyle
}%
\def\proofusing{
\eproofbit 
\setbox\proofrulename=\hbox\bgroup 
\let\wereinproofbit\proofusing
\kern0.3em$
}

\def\endprooftree{
\eproofbit 
  \dimen5 =0pt
\dimen0=\wd\proofabove \advance\dimen0-\shortenproofleft
\advance\dimen0-\shortenproofright
\dimen1=.5\dimen0 \advance\dimen1-.5\wd\proofbelow
\dimen4=\dimen1
\advance\dimen1\proofbelowshift \advance\dimen4-\proofbelowshift
\ifdim  \dimen1<0pt
\then   \advance\shortenproofleft\dimen1
        \advance\dimen0-\dimen1
        \dimen1=0pt
        \ifdim  \shortenproofleft<0pt
        \then   \setbox\proofabove=\hbox{%
                        \kern-\shortenproofleft\unhbox\proofabove}%
                \shortenproofleft=0pt
        \fi
\fi
\ifdim  \dimen4<0pt
\then   \advance\shortenproofright\dimen4
        \advance\dimen0-\dimen4
        \dimen4=0pt
\fi
\ifdim  \shortenproofright<\wd\proofrulename
\then   \shortenproofright=\wd\proofrulename
\fi
\dimen2=\shortenproofleft \advance\dimen2 by\dimen1
\dimen3=\shortenproofright\advance\dimen3 by\dimen4
\ifproofdots
\then
        \dimen6=\shortenproofleft \advance\dimen6 .5\dimen0
        \setbox1=\vbox to\proofdotseparation{\vss\hbox{$\cdot$}\vss}%
        \setbox0=\hbox{%
                \advance\dimen6-.5\wd1
                \kern\dimen6
                $\vcenter to\proofdotnumber\proofdotseparation
                        {\leaders\box1\vfill}$%
                \unhbox\proofrulename}%
\else   \dimen6=\fontdimen22\the\textfont2 
        \dimen7=\dimen6
        \advance\dimen6by.5\proofrulebreadth
        \advance\dimen7by-.5\proofrulebreadth
        \setbox0=\hbox{%
                \kern\shortenproofleft
                \ifdoubleproof
                \then   \hbox to\dimen0{%
                        $\mathsurround0pt\mathord=\mkern-6mu%
                        \cleaders\hbox{$\mkern-2mu=\mkern-2mu$}\hfill
                        \mkern-6mu\mathord=$}%
                \else   \vrule height\dimen6 depth-\dimen7 width\dimen0
                \fi
                \unhbox\proofrulename}%
        \ht0=\dimen6 \dp0=-\dimen7
\fi
\let\doll\relax
\ifwasinsideprooftree
\then   \let\VBOX\vbox
\else   \ifmmode\else$\let\doll=$\fi
        \let\VBOX\vcenter
\fi
\VBOX   {\baselineskip\proofrulebaseline \lineskip.2ex
        \expandafter\lineskiplimit\ifproofdots0ex\else-0.6ex\fi
        \hbox   spread\dimen5   {\hfi\unhbox\proofabove\hfi}%
        \hbox{\box0}%
        \hbox   {\kern\dimen2 \box\proofbelow}}\doll%
\global\dimen2=\dimen2
\global\dimen3=\dimen3
\egroup 
\ifonleftofproofrule
\then   \shortenproofleft=\dimen2
\fi
\shortenproofright=\dimen3
\onleftofproofrulefalse
\ifinsideprooftree
\then   \hskip.5em plus 1fil \penalty2
\fi
}

\newcounter{enums}

\newdimen\widelabel
\def\enumsentence{\@ifnextchar[{\@enumsentence}
{\refstepcounter{enums}\@enumsentence[(\theenums)]}}

\long\def\@enumsentence[#1]#2{\begin{list}{}{%
\advance\leftmargin by\widelabel \advance\labelwidth by \widelabel}
\item[#1] #2
\end{list}}

\newcounter{tempcnt}

\def\@item[#1]{\if@noparitem \@donoparitem
  \else \if@inlabel \indent \par \fi
         \ifhmode \unskip\unskip \par \fi 
         \if@newlist \if@nobreak \@nbitem \else
                        \addpenalty\@beginparpenalty
                        \addvspace\@topsep \addvspace{-\parskip}\fi
           \else \addpenalty\@itempenalty \addvspace\itemsep 
          \fi 
    \global\@inlabeltrue 
\fi
\everypar{\global\@minipagefalse\global\@newlistfalse 
          \if@inlabel\global\@inlabelfalse \hskip -\parindent \box\@labels
             \penalty\z@ \fi
          \everypar{}}\global\@nobreakfalse
\if@noitemarg \@noitemargfalse \if@nmbrlist \refstepcounter{\@listctr}\fi \fi
\setbox\@tempboxa\hbox{\makelabel{#1}}%
\global\setbox\@labels
 \hbox{\unhbox\@labels \hskip \itemindent
       \hskip -\labelwidth \hskip -\labelsep 
       \ifdim \wd\@tempboxa >\labelwidth 
                \box\@tempboxa
          \else \hbox to\labelwidth {\unhbox\@tempboxa}\fi
       \hskip \labelsep}\ignorespaces}

\newcounter{enumsi}

\newdimen\eeindent
\def\@mklab#1{\hfil#1}
\def\enummklab#1{\hfil(\eelabel)\hbox to \eeindent{\hfil#1}}
\def\enummakelabel#1{\enummklab{#1}\global\let\makelabel=\@mklab}
\def\toplabel#1{{\edef\@currentlabel{\p@enums\theenums}\label{#1}}}

\def\eenumsentence{\@ifnextchar[{\@eenumsentence}
{\refstepcounter{enums}\@eenumsentence[\theenums]}}

\long\def\@eenumsentence[#1]#2{\def\eelabel{#1}\let\holdlabel\makelabel%
\begin{list}{\alph{enumsi}.}{\usecounter{enumsi}%
\advance\leftmargin by \eeindent \advance\leftmargin by \widelabel%
\advance\labelwidth by \eeindent \advance\labelwidth by \widelabel%
\let\makelabel=\enummakelabel}
#2
\end{list}\let\makelabel\holdlabel}

\def\evnup{\@ifnextchar[{\@evnup}{\@evnup[0pt]}}

\def\@evnup[#1]#2{\setbox1=\hbox{#2}%
\dimen1=\ht1 \advance\dimen1 by -.5\baselineskip%
\advance\dimen1 by -#1%
\leavevmode\lower\dimen1\box1}

\makeatother

\begin{document}
\title{Intensional Models for the Theory of Types\thanks{The
Journal of Symbolic Logic, to appear.}}
\author{Reinhard Muskens}
\date{}
\maketitle
\begin{abstract}
In this paper we define \emph{intensional} models for the classical theory 
of types, thus arriving at an intensional type logic \textsf{ITL}.
Intensional models generalize Henkin's general models and have a natural
definition. As a class they do not validate the axiom of Extensionality. We
give a cut-free sequent calculus for type theory and show completeness of
this calculus with respect to the class of intensional models via a model
existence theorem. After this we turn our attention to applications.
Firstly, it is argued that, since \textsf{ITL} is truly intensional, it can
be used to model ascriptions of propositional attitude without predicting
\emph{logical omniscience}. In order to illustrate this a small fragment of
English is defined and provided with an \textsf{ITL} semantics. Secondly, it
is shown that \textsf{ITL} models contain certain objects that can be
identified with \emph{possible worlds}. Essential elements of modal logic
become available within classical type theory once the axiom of
Extensionality is given up.
\end{abstract}
\section{Introduction} 

The axiom scheme of Extensionality states that whenever two predicates or
relations are coextensive they must have the same properties:
\begin{equation}\label{ext}
\forall XY(\forall \vec{x}(X\vec{x}\leftrightarrow Y\vec{x})\to
\forall Z(ZX\to ZY))
\end{equation}
Historically Extensionality has always been problematic, the main problem
being that in many areas of application, though not perhaps in the
foundations of mathematics, the statement is simply false. This was
recognized by Whitehead and Russell in \emph{Principia Mathematica}
\cite{whitehead_russell}, where intensional functions such as `$A$ believes
that $p$' or `it is a strange coincidence that $p$' are discussed at length.
However, in the introduction to the second edition (1927) of the
\emph{Principia} Whitehead and Russell (influenced by Wittgenstein's
\emph{Tractatus}) already entertain the possibility that ``all functions of
functions are extensional''. Thirteen years later, in
Church's~\cite{chur:form40} canonical formulation of the Theory of Types, it
is observed that axioms of Extensionality should be adopted ``[i]n order to
obtain classical real number theory (analysis)'', a wording that does not
seem to rule out the option of not adopting them. Church's formulation of
type theory was completely syntactic and axioms could be adopted or dropped
at will, but in Henkin's~\cite{henk:comp50} classical proof of generalized
completeness the models that are considered, both the ``standard'' models
and the ``general'' ones, simply validate Extensionality. Although Henkin's
text still allows giving up the axiom\footnote{Henkin~\cite{henk:comp50}:
The axioms of extensionality \ldots\ can be dropped if we are willing to
admit models whose domains contain functions which are regarded as distinct
even though they have the same value for every argument.} the formal set-up
now effectively rules out intensional predicates and functions.

This poses problems for those areas of application of the logic where it is
important to distinguish between predicates that are coextensive and where
propositions that determine the same set of possible worlds should be kept
apart nevertheless. Linguistic semantics and Artificial Intelligence are
such applications and the problem has been dubbed one of ``logical
omniscience'' there, for it is with propositional attitudes like knowledge
and belief that predicates of predicates and predicates of propositions most
naturally arise. Is there a deep foundational difficulty with type theory
that makes the theory adequate for one area of application (mathematics) but
not for others? Or is it possible to come up with a revised and generalized
semantics for the logic, in which intensional predicates of predicates (or
intensional functions of functions) are allowed? In the latter case
Extensionality becomes a \emph{non}-logical axiom that can be added to the
theory for the purposes of one area of application while in other areas of
application it is not added.

Even if one is interested in mathematical applications of type theory only
there are good reasons to consider a generalization of its models in which
Extensionality fails. This was realized by Takahashi~\cite{taka:cut67} and
Prawitz~\cite{praw:haup68} in their (independent, but closely related)
proofs of Cut-elimination. These proofs make use of
what Andrews~\cite{andr:reso71} calls ``$V$-\emph{complexes}'', structures
whose typed domains consist of elements $\langle A,e\rangle$, where $A$ is a
term and $e$ is a possible extension of $A$. Clearly, two objects $\langle
A,e\rangle$ and $\langle B,e'\rangle$ can be distinct even if $e=e'$.
Andrews~\cite{andr:reso71} uses $V$-complexes to show that a certain
resolution system ${\mathcal R}$ corresponds to the first six axioms of
Church~\cite{chur:form40} (not comprising Extensionality).

$V$-complexes in themselves cannot be used as independent models for an
intensional type theory, as their definition depends on
Sch\"utte's~\cite{schu:synt60} ``semi-valuations'', essentially sets of
sentences (the ``$V$'' in ``$V$-complex'' ranges over semi-valuations). Is
it possible to define a stand-alone notion of general intensional model that
has $V$-complexes as a special case? I know of two proposals for such
general models, both recent. The first is found in
Fitting~\cite{fitt:type02}, the second in Benzm\"uller et al.~\cite{BBK04}.
In Fitting's ``generalized Henkin models'' abstraction may receive a
non-standard interpretation, while in the ``$\Sigma$-models'' of
Benzm\"uller et al.\ it is application that may be interpreted in a
non-standard way. Such non-standard evaluations seem unnecessary, however,
and in this paper, I will propose a simple definition of \emph{intensional
model} that generalizes Henkin models for type theory but gives all logical
operations their usual semantics. 
The system of type theory interpreted with the help
of these intensional models will be called \textsf{ITL}
(`Intensional Type Logic').

The rest of the paper is organized as follows. In the following section
we will consider some existing proposals to obtain intensionality and we
will argue that they all have a simple pattern in common that can be used to
obtain a general intensional logic. Section \ref{terms} gives the types and
terms of a type theory in the spirit of Church~\cite{chur:form40}
(but framed as a relational theory, as in Orey~\cite{orey:mode59} and
Sch\"utte~\cite{schu:synt60}). In section \ref{models} our notion of
intensional model will be defined, with a corresponding notion of
entailment. Section \ref{proofs} introduces a cut-free Gentzen calculus for
\textsf{ITL} while section \ref{modelex} proves a Model Existence theorem.
The proofs in that section all employ familiar techniques but are given as a
sanity check on the definition of the basic modeltheoretic notions. The last
two sections consider applications: Section \ref{ling} uses the logic to
provide a fragment of natural language with a truly intensional semantics
while section \ref{worlds} shows how possible worlds can be obtained as
certain objects in intensional models. A short conclusion ends the paper.

\section{Informal Analysis}\label{informal}
Given Leibniz's principle of the identity of indistinguishables and the
assumption that $\forall$ and $\to$ behave classically, $\forall Z(ZX\to
ZY)$ implies $X=Y$ and
\begin{equation}\label{ext2}
\forall XY(\forall \vec{x}(X\vec{x}\leftrightarrow Y\vec{x})\to
X=Y)
\end{equation}
therefore will be equivalent with Extensionality. This means that a semantics in
which this axiom fails cannot under reasonable assumptions identify the
semantic value of an expression with its extension, as $\forall
\vec{x}(X\vec{x}\leftrightarrow Y\vec{x})$ just states that $X$ and $Y$ are
co-extensive. 
The following is a propositional instantiation of (\ref{ext2}) (with the
length of $\vec{x}$ set to $0$ and $X$ and $Y$ instantiated as $\varphi$ and
$\psi$ respectively).
\begin{equation}\label{ext3}
(\varphi\leftrightarrow \psi)\to \varphi = \psi
\end{equation}
Here $\varphi\leftrightarrow \psi$ expresses that $\varphi$ and $\psi$ have
the same truth value whereas $\varphi = \psi$ says they are the same
proposition. Typically we want this scheme to fail, as sentences with the
same truth values may be distinguishable in the sense that one is believed
while the other is not, or that one is a strange coincidence, while the
other is entirely expected, etc. Many propositions must therefore be allowed
to exist, although each proposition must assume one of two truth values if
we want to retain classicality.

Although the semantic values of sentences cannot be \emph{equated} with
their truth values, it still seems reasonable to require that they should
\emph{determine} these truth values, while the values of expressions of
higher type should likewise determine their extensions. If this is accepted
the picture that arises is that logical expressions are sent to their
(intensional) values by some function $I$, while these values in their turn
are connected with extensions by a function $E$. The latter typically does
not need to be injective. Readers familiar with the extensive literature on
intensionality (sometimes dubbed \emph{hyper}intensionality) will be aware
that, while there are many divergent proposals for what intensions
\emph{are}, the pattern just sketched is well-nigh ubiquitous. Already in
Frege's~\cite{freg:sinn92} pivotal work an expression \emph{expresses} a
sense (the function $I$) while a sense in its turn \emph{determines} a
reference ($E$). In modal logic, intensions are functions from possible
worlds to extensions (let us call such functions \emph{modal intensions})
and the function $E$ can be viewed as application of such functions to a
fixed `actual' world $w_0$. It is obvious that this modal strategy does not
individuate intensions finely enough (essentially since, if $W_1$ and $W_2$
denote sets of possible worlds or the characteristic functions of such sets,
$\forall w(W_1 w\leftrightarrow W_2 w)\to W_1=W_2$ will hold) and many
researchers have sought notions of intension coming with more fine-grained
criteria of individuation than modal intensions come with.
Carnap~\cite{carn:mean}, for example, who defined a precursor to the now
usual possible worlds analysis of natural language, noticed the
problems that this analysis suffers from and proposed a theory of
\emph{structured meanings} that was later worked out in
Lewis~\cite{lewi:gene72} and Cresswell~\cite{cres:struc85}. On Lewis'
account the \emph{meaning} of an expression is a finite ordered tree having
at each node a category and an appropriate modal intension. The modal
intension at the root node of this tree is the one associated with the
expression as a whole and so in this theory the function $I$ assigns finite
ordered trees of categories and modal intensions to any given expression,
while $E$ is the function that takes any such tree and returns $W w_0$,
where $W$ is the modal intension found at its root node. 

Another example of an approach in which the functions $I$ and $E$ can easily
be recognised is the theory of \emph{impossible worlds}. The idea behind
this line of thought is that if the usual set of possible worlds is not
large enough to make enough distinctions between semantic values, extra
worlds, impossible ones, should be added. A key point is that the logical
operators need \emph{not} have their usual meaning at these points of
reference and that logical validities will therefore cease to hold
throughout the set of all worlds. The name ``impossible (possible) world''
derives from Hintikka~\cite{hint:impo75}, but the idea was also present in
Montague~\cite{mont:univ70} and Cresswell~\cite{cres:inte72} and has been followed
up in Rantala~\cite{rant:quan82}, Barwise~\cite{barw:info97}, and 
Zalta~\cite{zalt:clas97}, for example. The function $I$ here is the
function that sends each sentence to a set of possible and impossible
worlds, wheras $E$ can be described as $\lambda W.Ww_0$, as in the ordinary
modal account.

Other approaches to intensionality may have different conceptions of the
nature of intensions, but will also follow a two-stage pattern in which
expressions are first sent to their intensions (whatever these are) and
intensions are subsequently related to extensions. For example, in Property
Theory (Turner~\cite{turn:theo87}, Chierchia and Turner~\cite{chie:sema88})
one finds a homomorphism $T$ sending an algebra of `information units'
$\mathbf{I}$ to a boolean algebra $\mathbf{P}$ (Chierchia and Turner
\cite{chie:sema88}). Here the information units act as intensions while the
elements of $\mathbf{P}$ are the extensions and $T$ is the function we have
called $E$. Thomason~\cite{thom:mode80} uses a higher-order logic that obeys
a form of Extensionality to interpret natural language
sentences in a domain $D_p$ of propositions considered as primitive entities
and then uses a function $^\cup$ (our $E$) to send these propositions to
their extensions (see also Muskens~\cite{musk:sense05}). Moschovakis
\cite{mosc:sens94}, to give a last example, identifies senses with
algorithms and references with the values that these algorithms return. Here
the function that sends expressions to algorithms is our $I$ while $E$
assigns to each algorithm the value returned.

Thus while opinions about the nature of intensions radically diverge, all
proposals follow a simple two-stage pattern. The aim of this paper is not to
add one more theory of intension to the proposals that have already been
made, but is an investigation of their common underlying logic. The idea
will be that the two-stage set-up is essentially \emph{all} that is needed
to obtain intensionality. For the purposes of logic it suffices to consider
intensions as abstract objects; the question what intensions \emph{are},
while philosophically important, can be abstracted from. Conversely, while
many positions regarding the ultimate nature of intensions seem rationally
possible and no knock-down arguments are likely to decide the matter, it
does equally seem possible to rationally converge on a logic describing what
intensions \emph{do}. Here we attempt to contribute to that logic.

\section{Terms}\label{terms}
In this section the types and terms of \textsf{ITL} will be defined and some
notation will be adopted. While this \textsf{ITL} syntax will be given an
intensional interpretation in the next section, it essentially is the syntax
of the simple type theory of Church's~\cite{chur:form40}. The intended 
interpretation will be \emph{relational}, however, as in
Orey~\cite{orey:mode59} and Sch\"utte~\cite{schu:synt60}, not
\emph{functional}, as in Church's original work.

Assuming that some finite set ${\mathcal B}$ of basic types is given, the
following definition gives the set of types.
\begin{definition}
The set ${\mathcal T}$ of \emph{types} is the smallest set of strings
such that
\begin{enumerate}
\item ${\mathcal B}\subseteq{\mathcal T}$
\item If $\alpha_1,\ldots,\alpha_n\in{\mathcal T}$ ($n\geq 0$) then
$\langle\alpha_1\ldots\alpha_n\rangle\in{\mathcal T}$
\end{enumerate}
\end{definition}
\noindent Types formed with the second clause of this definition will be
called \emph{complex}. The intended interpretation is that the extension of
an object of type $\langle\alpha_1\ldots\alpha_n\rangle$ is an $n$-place
relation taking objects of type $\alpha_i$ as its $i$-th argument. Note
that, as a limiting case, $\langle\rangle$ is defined to be a (complex)
type; this will be the type of \emph{propositions}, with \emph{truth values}
as extensions.

A \emph{language} will be a countable set of uniquely typed non-logical
constants. If ${\mathcal L}$ is a language, the set of constants from
${\mathcal L}$ that have type $\alpha$ is denoted ${\mathcal L}_\alpha$. For
each $\alpha\in{\mathcal T}$ we moreover assume the existence of a
denumerably infinite set ${\mathcal V}_\alpha$ of variables with unique type
$\alpha$. We let ${\mathcal V}=\bigcup_{\alpha\in{\mathcal T}}{\mathcal
V}_\alpha$. 

The following definition gives us terms in all types. Apart from variables
and non-logical vocabulary there will be a sentence $\bot$ that is always
false, and there will be application and abstraction. Furthermore, a symbol 
$\subset$ will denote \emph{inclusion of extensions}, so that 
$A\subset B$ is true if the extension of $A$ is a subset of that of $B$.

\begin{definition}
Let ${\mathcal L}$ be a language. Define sets $T_\alpha^{\mathcal L}$ of 
\emph{terms} of ${\mathcal L}$ of type $\alpha$, for each
$\alpha\in{\mathcal T}$, as follows.
\begin{enumerate}
\item ${\mathcal L}_\alpha\subseteq T_\alpha^{\mathcal L}$ and
${\mathcal V}_\alpha\subseteq T_\alpha^{\mathcal L}$ for
each $\alpha\in{\mathcal T}$
\item $\bot\in T_{\langle\rangle}^{\mathcal L}$
\item If $A\in T_{\langle\alpha_1\alpha_2\ldots\alpha_n\rangle}^{\mathcal L}$ 
and $B\in T_{\alpha_1}^{\mathcal L}$, then 
$(AB)\in T_{\langle\alpha_2\ldots\alpha_n\rangle}^{\mathcal L}$
\item If $A\in T_{\langle\alpha_2\ldots\alpha_n\rangle}^{\mathcal L}$
and $x\in {\mathcal V}_{\alpha_1}$, then
$(\lambda x.A)\in T_{\langle\alpha_1\alpha_2\ldots\alpha_n\rangle}^{\mathcal L}$
\item If $A\in T_\alpha^{\mathcal L}$ and
$B\in T_\alpha^{\mathcal L}$ then 
$(A\subset B)\in T_{\langle\rangle}^{\mathcal L}$, if $\alpha$ is complex
\end{enumerate}
\end{definition}
\noindent We will write $T^{\mathcal L}$ for the set of \emph{all} terms of the
language ${\mathcal L}$, i.e.\ for the union $\bigcup_{\alpha\in{\mathcal
T}}T_\alpha^{\mathcal L}$. If $A$ is a term of type $\alpha$, we may
indicate this by writing $A_\alpha$ and we will use $\varphi$, $\psi$,
$\chi$ for terms of type $\langle\rangle$, which we call \emph{formulas}.
The notions \emph{free} and \emph{bound} occurrence of a variable and the
notion $B$ \emph{is free for} $x$ \emph{in} $A$ are defined as usual, as are
\emph{closed} terms and \emph{sentences}. \emph{Substitutions} are functions
$\sigma$ from variables to terms such that $\sigma(x)$ has the same type as
$x$. If $\sigma$ is a substitution then the substitution $\sigma'$ such that
$\sigma'(x)=A$ and $\sigma'(y)=\sigma(y)$ for all $y\not\equiv x$ is denoted
as $\sigma[x:=A]$. If $A$ is a term and $\sigma$ is a substitution,
$A\sigma$, the extension of $\sigma$ to $A$, is defined in the usual way.
The substitution $\sigma$ such that $\sigma(x_i)=A_i$ and $\sigma(y)=y$ if
$y\notin\{x_1,\ldots,x_n\}$ is written as $\{x_1:=A_1,\ldots,x_n:=A_n\}$.
Parentheses in terms will often be dropped on the understanding that $ABC$
is $((AB)C)$, i.e.\ association is to the left.

Our stock of operators may seem somewhat spartan, but is rich enough to let
the usual connectives and quantifiers be defined. In particular, $\forall$,
$\to$ and $=$ are easily obtained.

\begin{definition}\label{abbr}
Write
\vspace{-2ex}
\begin{align*}
\varphi\to\psi
&\quad\mbox{for}\quad\varphi\subset\psi,\\
\top&\quad\mbox{for}\quad\bot\to\bot,\\
\forall x\varphi&\quad\mbox{for}\quad
(\lambda x.\top)\subset(\lambda x.\varphi)\mbox{, and}\\
A_\alpha=B_\alpha&\quad\mbox{for}\quad
\forall x_{\langle\alpha\rangle}\,(xA\to xB).
\end{align*}
\end{definition}
\noindent The operators $\lnot$, $\land$, $\lor$, $\leftrightarrow$ and 
$\exists$ are defined as usual. 

Our presentation of the logic will revolve around sequents.
A \emph{signed sentence} of ${\mathcal L}$ will be a pair
$\langle \mathbf{\sf L},\varphi\rangle$ (written $\T\varphi$) or a pair
$\langle \mathbf{\sf R},\varphi\rangle$ (written $\F\varphi$), such
that $\varphi$ is a sentence of ${\mathcal L}$ ($\mathbf{\sf L}$ indicates
`left' and $\mathbf{\sf R}$ indicates `right').
A \emph{sequent}
of ${\mathcal L}$ is a set of signed 
sentences of ${\mathcal L}$. Letting
sequents be sets has some advantages, but we may also want to use a more
conventional form and write 
$\Pi\Rightarrow\Sigma$ for
$\{\T\varphi\mid\varphi\in\Pi\}\cup\{\F\varphi\mid\varphi\in\Sigma\}$ if
$\Pi$ and $\Sigma$ are sets of sentences.

\section{Intensional Models}\label{models}

Let us turn to the semantics of \textsf{ITL}, which will essentially follow
the two-stage pattern discussed above. The following definition sets up the
usual hierarchies of objects and provides some of the usual notation.

\begin{definition}\label{collection}
A \emph{collection of domains}
will be a set $\{D_\alpha\mid\alpha\in{\mathcal T}\}$ of pairwise disjoint
non-empty sets. An \emph{assignment} $a$ for a collection of domains 
$D=\{D_\alpha\mid\alpha\in{\mathcal T}\}$ is a function which has the set of
variables ${\mathcal V}$ as domain and has the property that $a(x)\in
D_\alpha$ if $x\in{\mathcal V}_\alpha$. The set of all assignments for $D$
is denoted ${\mathcal A}_D$. If $a$ is an assignment, $d\in D_{\alpha}$, and
$x$ is a variable of type $\alpha$, $a[d/x]$ is defined by letting
$a[d/x](x)=d$ and $a[d/x](y)=a(y)$, if $y$ is not equal to $x$.
\end{definition}
\noindent Note that we have not imposed any non-trivial relations between
the elements of any given collection of domains $D$. In particular
we have \emph{not} required domains
$D_{\langle\alpha_1\ldots\alpha_n\rangle}$ to consist of relations over
lower domains. This is because we need to tease apart the intensions and
extensions of terms of complex type. While extensions of such terms will be
certain relations, with their identity criteria therefore given by set
membership, the \emph{intension functions} defined below send terms to
almost arbitrary domain elements. 

\begin{definition}\label{intensionfunction}
An \emph{intension function} for a collection of domains
$D=\{D_\alpha\mid\alpha\in{\mathcal T}\}$ and a language
${\mathcal L}$ is a function 
$I\colon{\mathcal A}_D\times T^{\mathcal L}\to D$ such that
\begin{enumerate}
\item $I(a,A)\in D_\alpha$, if $A$ is of type $\alpha$
\item $I(a,x)=a(x)$, if $x$ is a variable
\item $I(a, A)=I(a', A)$, if $a$ and $a'$ agree on all variables free in $A$
\item $I(a,A\{x:=B\})=I(a[I(a,B)/x],A)$, if $B$ is free for
$x$ in $A$
\end{enumerate}
\end{definition}
\noindent Intension functions are the formal counterpart of the functions
$I$ that were discussed informally above. They take an extra assignment
argument in order to take care of free variables.

The second part of our formalisation of the two-stage architecture discussed
above are the \emph{extension} functions of definition
\ref{extensionfunctions}. They send objects of complex types to certain
relations over the relevant domains.  We first give very general constraints
and will put more requirements on useful extension functions in definition
\ref{genmod}.
\begin{definition}\label{extensionfunctions}
An \emph{extension function} for a collection of domains
$D=\{D_\alpha\mid\alpha\in{\mathcal T}\}$ is a function $E$
with domain $\cup\{D_\alpha\mid\alpha\mbox{ is complex}\}$
such that $E(d)\subseteq D_{\alpha_1}\times\cdots\times D_{\alpha_n}$
whenever $d\in D_{\langle\alpha_1\ldots\alpha_n\rangle}$. 
\end{definition}
\noindent The restriction of $E$ to $D_\alpha$ is written as
$E_\alpha$, for any complex type $\alpha$, so that 
$E_\alpha\colon D_\alpha\to{\mathcal P}(D_{\alpha_1}\times\cdots\times
D_{\alpha_n})$ if $\alpha=\langle\alpha_1\ldots\alpha_n\rangle$.

The limiting case that $n=0$ is of some interest here. In this case the product
$D_{\alpha_1}\times\cdots\times D_{\alpha_n}$ equals $\{\langle\rangle\}$.
We identify $\langle\rangle$ with $\varnothing$, $\varnothing$ with $0$, and
$\{\varnothing\}$ with $1$, so that 
$E_{\langle\rangle}\colon D_{\langle\rangle}\to\{0,1\}$
if $E_{\langle\rangle}$ is an extension function of type
$\langle\rangle$ for $D$. Note that, while the range of $E_{\langle\rangle}$
thus consists of the two standard truth-values, the domain
$D_{\langle\rangle}$ of propositions can have any cardinality $\geq 2$.
Propositions with the same truth-value need not be identified and, as will
become apparent, even propositions that receive the same truth value
in all structures need not be identical in any given structure.
\begin{definition}
A \emph{generalized frame} for the language ${\mathcal L}$
is a triple $\langle D,I,E\rangle$ such
that $D$ is a collection of domains, $I$ is an intension
function for $D$ and ${\mathcal L}$, and $E$ is an extension
function for $D$. 
\end{definition}
\noindent We are interested in the extensions $E(I(a,A_\alpha))$ of
terms $A$ of complex type $\alpha$. Let $V$ be the composition of
$E$ and $I$, so that, in the interest of readability, we can write
$V(a,A)$, for $E(I(a,A))$. The following
definition, which gives the central notion of this paper, puts constraints
on intension and extension functions that cause terms to get their usual
semantic values.
\begin{definition}\label{genmod}
A generalized frame $\langle D,I,E\rangle$
for ${\mathcal L}$ is an \emph{intensional model} for ${\mathcal L}$ if
\begin{enumerate}
\item $V(a,\bot)=0$
\item $V(a,AB)=\{\langle \vec{d}\rangle\mid
\langle I(a,B),\vec{d}\rangle\in V(a,A)\}$
\item $V(a,\lambda x_\beta.A)=\{\langle d,\vec{d}\rangle\mid
d\in D_\beta\mbox{ and }\langle \vec{d}\rangle\in V(a[d/x],A)\}$
\item $V(a,A\subset B)=1\Longleftrightarrow 
V(a,A)\subseteq V(a,B)$
\end{enumerate}
\end{definition}
\noindent To better understand the motivation behind the second and third
clauses of this definition, it may help to consider that any $n+1$ place
relation $R$ can be thought of as a unary function $F$ such that
$F(d)=\{\langle\vec{d}\rangle\mid\langle d,\vec{d}\rangle\in R\}$. Thus 
$V(a,AB)=F(I(a,B))$, where $F$ is the function corresponding to $V(a,A)$ and
$V(a,\lambda x_\beta.A)$ corresponds to the function $F$ such that $F(d)=
V(a[d/x],A)$ for each $d\in D_\beta$. For further discussion of this little
trick in an extensional setting see Muskens~\cite{musk:rela89,mp}.

If $M=\langle D,I,E\rangle$ is an intensional model, $a$ is an assignment for
$D$, and $\varphi$ is a formula, we may alternatively write
$M\models\varphi[a]$ for $V(a,\varphi)=1$.  In case $\varphi$ is a
sentence it makes sense to write $M\models\varphi$ if $M\models\varphi[a]$
for some $a$. The following facts are unsurprising but useful.

\begin{proposition} \label{facts}
Let $M=\langle D,I,E\rangle$ be an intensional model, 
and let $a$ be an assignment for $D$. Then, for all $\varphi$, $\psi$, $A$,
$B$ and $B'$ of appropriate types,
\begin{enumerate}
\item $V(a,\varphi\to\psi)=0$ iff $V(a,\varphi)=1$ and
$V(a,\psi)=0$;
\item $V(a,\forall x_\alpha\varphi)=1$ iff
$V(a[d/x],\varphi)=1$ for all $d\in D_\alpha$;
\item $V(a,(\lambda x.A)B)=V(a,A\{x:=B\})$, if $B$ is
free for $x$ in $A$;
\item If $V(a,A=B)=1$ then $V(a,A\subset B)=1$;
\item $V(a,A=A)=1$;
\item If $V(a,B=B')=1$ then $V(a,A\{x:=B\}=A\{x:=B'\})=1$, provided
$B$ and $B'$ are free for $x$ in $A$.
\end{enumerate}
\end{proposition}
\begin{proof} Left to the reader.\end{proof}
\noindent Note that $\beta$-conversion preserves \emph{extensional}
identity, but that it does not necessarily preserve \emph{intensional}
identity, i.e.\ $(\lambda x_\alpha.A)B=A\{x:=B\}$ is not necessarily true
given the usual side condition. Similar remarks can be made about
$\eta$-conversion and even about $\alpha$-conversion. Since it is not
necessary to hardwire these principles into the logic, we have chosen not to
do so. However, the principles can clearly be added to the logic by means of
an axiomatic extension. In section \ref{ling} below, where a linguistic
application is considered, this axiomatic extension will be given.

The last two statements in proposition \ref{facts} above show that $=$ is
the usual congruence, but intensional models may still have the undesirable
property that $=$ does not denote true identity of intension. This is an
anomaly we want to get rid of. Intensional models are called \emph{normal}
just in case they have the desired property.
\begin{definition}
An intensional model $M=\langle D,I,E\rangle$ is \emph{normal} if, for any
type $\alpha$, any $d,d'\in D_\alpha$, and any $a$,
$\langle d,d'\rangle\in V(a,\lambda x_\alpha\lambda x'_\alpha.x=x')$
implies $d=d'$.
\end{definition}
\noindent That a restriction to normal intensional models does not buy us
any new truths is shown by the next proposition. Its proof uses the Axiom of
Choice unless $M$ is countable.
\begin{proposition}\label{normalize}
Let $M$ be an intensional model. There is a normal intensional model
$\overline{M}$ such that $M\models\varphi\Longleftrightarrow \overline{M}\models\varphi$
for each sentence $\varphi$.
\end{proposition}
\begin{proof}
Suppose $M=\langle D,I,E\rangle$. Let $\sim$ be the relation given
by $d\sim d'$ iff $\langle d,d'\rangle\in 
V(a,\lambda x_\alpha\lambda x'_\alpha.x=x')$ for any 
$d,d'\in D_\alpha$ and any $\alpha$ (where $a$ is arbitrary). 
Clearly, $\sim$ is an equivalence
relation. Note that, by proposition \ref{facts} and
definition \ref{intensionfunction}, for any term $A$,
\begin{equation}\label{congr}
d\sim d' \Longrightarrow I(a[d/x],A)\sim I(a[d'/x],A)\ .
\end{equation}
Define $\overline{d}=\{d'\mid d\sim d'\}$, let 
$\overline{D}_\alpha=\{\overline{d}\mid d\in D_\alpha\}$, and let
$\overline{D}=\{\overline{D}_\alpha\mid\alpha\in{\mathcal T}\}$. Let $f$
be a function such that $f(\overline{d})\in\overline{d}$, if
$\overline{d}\in \overline{D}_\alpha$. For
any assignment $a$ for $\overline{D}$, let $a^\circ$ be the assignment
for $D$ defined by $a^\circ(x)=f(a(x))$, for all $x$. Let
$\overline{I}(a,A)=\overline{I(a^\circ,A)}$, for each assignment $a$ for
$\overline{D}$ and each term $A$. Then $\overline{I}$ is an intension
function for $\overline{D}$. The first three requirements of definition
\ref{intensionfunction} are easily checked, so let us check the last
requirement. Note that
\begin{eqnarray*}
I(a^\circ,A\{x:=B\})&=&\mbox{(Definition \ref{intensionfunction})}\\
I(a^\circ[I(a^\circ,B)/x],A)&\sim&\mbox{(\ref{congr})}\\
I(a^\circ[f(\overline{I(a^\circ,B)})/x]),A)
&=&\mbox{(definition of $\overline{I}$)}\\
I(a^\circ[f(\overline{I}(a,B))/x]),A)
&=&\mbox{(definition of $\circ$)}\\
I((a[\overline{I}(a,B)/x])^\circ,A)\ .
\end{eqnarray*}
From this conclude that $\overline{I}(a,A\{x:=B\})=\overline{I}(a[\overline{I}(a,B)/x],A)$.

Define $\overline{E}$ by
letting
$\overline{E}(\overline{d}_\alpha)=
\{\langle\overline{d_1},\ldots,\overline{d_n}\rangle\mid
\langle d_1,\ldots, d_n\rangle\in E(d)\}$,
if $\alpha$ is complex. It is easy to see that this is
well-defined. Since $\langle\overline{d_1},\ldots,\overline{d_n}\rangle\in
\overline{E}(\overline{I}(a,A))$ iff $\langle d_1,\ldots, d_n\rangle\in
V(a^\circ,A)$ it follows that
$\overline{M}=\langle\overline{D},\overline{I},\overline{E}\rangle$
is an intensional model, $M\models\varphi\Longleftrightarrow \overline{M}\models\varphi$
for each sentence $\varphi$, and $\overline{M}$ is normal.
\end{proof}

\noindent Now that the situation with respect to normality and non-normality
of intensional models has become clear, we can define our semantic notion of
consequence.

\begin{definition}
An intensional model $M$ for ${\mathcal L}$ \emph{refutes} a sequent
$\Pi\Rightarrow\Sigma$ of ${\mathcal L}$ if $M\models\varphi$ for all
$\varphi\in\Pi$ and $M\not\models\varphi$ for all $\varphi\in\Sigma$. A
sequent $\Gamma$ is \emph{i-valid} if no
intensional model for ${\mathcal L}$ refutes $\Gamma$. $\Pi$ \emph{i-entails}
$\Sigma$, $\Pi\models_i\Sigma$, if $\Pi\Rightarrow\Sigma$ is i-valid.
\end{definition}
\noindent Let us take stock. We have defined a notion of intensional model
following the two-stage pattern discussed in section \ref{informal}. This is
also the pattern followed in Fitting~\cite{fitt:type02}, but we have avoided
the complex ``abstraction designation functions'' that are used there but
do not seem to have a justification beyond the fact that they are
needed in proofs. Intensional models are a further generalisation of Henkin
models in the following sense. While in intensional models the functions
\[E_{\langle\alpha_1\ldots\alpha_n\rangle}\colon 
D_{\langle\alpha_1\ldots\alpha_n\rangle}\to{\mathcal
P}(D_{\alpha_1}\times\cdots\times D_{\alpha_n})\] need neither be injective
nor surjective, the usual Henkin models are essentially obtained if an
injectivity requirement is imposed. An additional requirement of
surjectivity brings us to a variant of the so-called \emph{standard} models
of type theory.

\begin{table}[t]
\begin{center}$
\begin{array}{cc}\hline \\
\multicolumn{2}{c}{
\prooftree
\Pi\Rightarrow\Sigma\justifies
\Pi'\Rightarrow\Sigma'
\using [W]
\endprooftree,\quad\mbox{ if }\Pi\subseteq\Pi',\ \Sigma\subseteq \Sigma'
}\\[4ex]
\prooftree
\justifies
\Pi,\varphi\Rightarrow\Sigma,\varphi
\using [R]
\endprooftree&
\prooftree
\justifies
\Pi,\bot\Rightarrow\Sigma
\using [\bot \mathbf{\sf L}]
\endprooftree\\[4ex]
\prooftree
\Pi,A\{x:=B\}\vec{C}\Rightarrow\Sigma\justifies
\Pi,(\lambda x.A)B \vec{C}\Rightarrow\Sigma
\using [\lambda \mathbf{\sf L}]
\endprooftree&
\prooftree
\Pi\Rightarrow \Sigma,A\{x:=B\}\vec{C}\justifies
\Pi\Rightarrow\Sigma,(\lambda x.A)B \vec{C}
\using [\lambda \mathbf{\sf R}]
\endprooftree\\[3ex]
\mbox{if $B$ is free for $x$ in $A$}&
\mbox{if $B$ is free for $x$ in $A$}\\[4ex]
\prooftree
\Pi,B\vec{C}\Rightarrow \Sigma\qquad
\Pi\Rightarrow\Sigma,A\vec{C}\justifies
\Pi,A\subset B\Rightarrow\Sigma
\using [\mathop{\subset} \mathbf{\sf L}]
\endprooftree&
\prooftree
\Pi,A\vec{c}\Rightarrow \Sigma,B\vec{c}\justifies
\Pi\Rightarrow\Sigma,A\subset B
\using [\mathop{\subset} \mathbf{\sf R}]
\endprooftree\\[3ex]
&\mbox{if the constants $\vec{c}$ are fresh}\\ \hline
\end{array}
$\end{center}
\caption{\label{gentzen}Gentzen rules for \textsf{ITL}.}
\end{table}

\section{Proof Theory}\label{proofs}
We now provide the relation of i-entailment with what
will turn out to be a syntactic equivalent. The rules in Table
\ref{gentzen}, for which the usual notational conventions apply, constitute
a Gentzen sequent calculus for \textsf{ITL}. If $\Pi\Rightarrow\Sigma$ is a
(finite or infinite) sequent, then we say that $\Pi\Rightarrow\Sigma$ is
\emph{provable}, $\Pi\vdash\Sigma$, if there are finite $\Pi_0\subseteq\Pi$
and $\Sigma_0\subseteq\Sigma$ such that $\Pi_0\Rightarrow\Sigma_0$ can be
proved in this calculus. The following theorem states that the calculus is
sound.

\begin{theorem}[Soundness]
If a sequent $\Gamma$ is provable, $\Gamma$ is i-valid. Hence
$\Pi\vdash\Sigma\Longrightarrow\Pi\models_i\Sigma$
\end{theorem}
\begin{proof}Left to the reader. (The proof involves some observations about
the behaviour of intension functions when the language is extended.)\end{proof}

\begin{table}[t]
\begin{center}$
\begin{array}{cc}\hline \\[-1ex]
\prooftree
\justifies
\Pi\Rightarrow\Sigma,\top
\using [\top \mathbf{\sf R}]
\endprooftree\\[4ex]
\prooftree
\Pi,\psi\Rightarrow \Sigma\qquad
\Pi\Rightarrow\Sigma,\varphi\justifies
\Pi,\varphi\to\psi\Rightarrow\Sigma
\using [\mathop{\to} \mathbf{\sf L}]
\endprooftree&
\prooftree
\Pi,\varphi\Rightarrow \Sigma,\psi\justifies
\Pi\Rightarrow\Sigma,\varphi\to\psi
\using [\mathop{\to} \mathbf{\sf R}]
\endprooftree\\[4ex]
\prooftree
\Pi,\varphi\{x:=A\}\Rightarrow\Sigma\justifies
\Pi,\forall x\varphi\Rightarrow\Sigma
\using [\forall \mathbf{\sf L}]
\endprooftree&
\prooftree
\Pi\Rightarrow \Sigma,\varphi\{x:=c\}\justifies
\Pi\Rightarrow\Sigma, \forall x\varphi
\using [\forall \mathbf{\sf R}]
\endprooftree\\[3ex]
&\mbox{where $c$ is fresh}\\[2ex]

\prooftree
\Pi,A\doteq B\Rightarrow \Sigma, \varphi\{x:=A\}\justifies
\Pi,A\doteq B\Rightarrow \Sigma, \varphi\{x:=B\}
\using [\mathop{=} \mathbf{\sf L}]
\endprooftree&
\prooftree
\justifies
\Pi\Rightarrow\Sigma,A=A
\using [\mathop{=} \mathbf{\sf R}]
\endprooftree\\[3ex]
\mbox{where $A\doteq B$ is $A=B$ or $B=A$}\\ \hline
\end{array}
$\end{center}
\caption{\label{derived}Some classical rules derivable in \textsf{ITL}.}
\end{table}

\noindent That the converse (generalized completeness) also holds will be
shown in the next section. 

While the rules in Table \ref{gentzen} suffice to characterize the
$\models_i$ relation, it is pleasant to also have the usual classical
Gentzen rules for the defined connectives at one's disposal. These are
available as derived rules. By way of example those for $\top$, $\to$,
$\forall$, and $=$ are given in Table \ref{derived}. Given the abbreviations
in definition \ref{abbr}, they are easily derivable from the \textsf{ITL}
rules, as the reader may verify. Note that in view of the correctness of
these rules it seems reasonable to say that \textsf{ITL} is indeed a
\emph{classical} logic. 
%
%
\section{Model Existence}\label{modelex}
The purpose of this section---which could be skipped on a first reading by
readers mainly interested in the general characteristics of our logic---is
to prove Generalized Completeness and some of its friends, such as the
generalized L\"owenheim-Skolem and Compactness theorems. We will do this in
the way Smullyan~\cite{smul:firs68} did it for first-order logic, via a
central \emph{Model Existence} theorem from which the desired theorems all
follow as corollaries. First it will be proved that certain ``Hintikka''
sequents, which can be thought of as resulting from a systematic but
unsuccessful attempt to construct a Gentzen proof from the bottom up, are
refutable. This is then used to show refutability of a wide class of
sequents.

The definition of Hintikka sequents is close to that of the ``Hintikka
sets'' in Smullyan~\cite{smul:firs68} and
Fitting~\cite{fitting96,fitt:type02}, but is also analogous to that of
Sch\"utte's~\cite{schu:synt60} semi-valuations.
\begin{definition}\label{hintikka}
A sequent $\Gamma$ of ${\mathcal L}$ is called a \emph{Hintikka sequent} 
in ${\mathcal L}$ if the following hold:
\begin{enumerate}
\item $\{\T\varphi,\ \F\varphi\}\not\subseteq\Gamma$ for any
sentence $\varphi$;
\item $\T\bot\notin\Gamma$;
\item $\T (\lambda x.A)B\vec{C}\in\Gamma\Longrightarrow
\T A\{x:= B\}\vec{C}\in\Gamma$, if
$\lambda x.A$,
$B$, and the sequence of terms $\vec{C}$ are closed and of appropriate
type;
\item $\F (\lambda x.A)B\vec{C}\in\Gamma\Longrightarrow
\F A\{x:= B\}\vec{C}\in\Gamma$, if
$\lambda x.A$,
$B$, and the sequence of terms $\vec{C}$ are closed and of appropriate
type;
\item $\T A\subset B\in\Gamma\Longrightarrow\T B\vec{C}\in\Gamma$ or 
$\F A\vec{C}\in\Gamma$, for all closed $A$, $B$ and sequences of
closed $\vec{C}$ of appropriate types;
\item $\F A\subset B\in\Gamma\Longrightarrow{}$ there are constants
$\vec{c}$ of appropriate types such that $\{\T A\vec{c},\F B\vec{c}\}
\subseteq\Gamma$.
\end{enumerate}
A Hintikka sequent $\Gamma$ in ${\mathcal L}$ is said to be \emph{complete} if
$\T\varphi\in\Gamma$ or
$\F\varphi\in\Gamma$, for each sentence $\varphi$ of ${\mathcal L}$.
\end{definition}
\noindent A key property of Hintikka sequents is that they are refuted by
intensional models, as the following lemma shows. The intensional model
constructed in its proof is closely akin to Andrews' $V$-complexes.

\begin{lemma}[Hintikka Lemma]
Let $\Gamma$ be a Hintikka sequent in a language ${\mathcal L}$ such that
${\mathcal L}_\alpha\neq\varnothing$ if $\alpha$ is basic. Then $\Gamma$ is
refuted by an intensional model. If $\Gamma$ is complete, then $\Gamma$ is
refuted by a normal countable intensional model.
\end{lemma}
\begin{proof}

Let $\Gamma$ be a Hintikka sequent in the language ${\mathcal L}$ as
described.  We will find an intensional model refuting $\Gamma$ using the
Takahashi-Prawitz construction. The following induction on type complexity
defines domains $D_\alpha$ as sets of pairs $\langle A,e\rangle$, where $A$
is a closed term of type $\alpha$ and $e$ is called a \emph{possible
extension} of $A$.
\begin{enumerate}
\item If $\alpha$ is basic let
$D_\alpha=\{\langle c,c\rangle\mid c\in{\mathcal L}_\alpha\}$; 
\item If $\alpha=\langle\alpha_1\ldots\alpha_n\rangle$ let
$\langle A_\alpha,e\rangle\in D_\alpha$ iff $A$ is closed,
$e\subseteq D_{\alpha_1}\times\cdots\times D_{\alpha_n}$
and, whenever $\langle B_1,e_1\rangle\in D_{\alpha_1},\ldots,
\langle B_n,e_n\rangle\in D_{\alpha_n}$
\begin{enumerate}
\item If $\T AB_1\ldots B_n\in\Gamma$ then 
$\langle \langle B_1,e_1\rangle,\ldots,\langle B_n,e_n\rangle\rangle\in e$;
\item If $\F AB_1\ldots B_n\in\Gamma$ then 
$\langle \langle B_1,e_1\rangle,\ldots,\langle B_n,e_n\rangle\rangle\notin e$.
\end{enumerate}
\end{enumerate}
It is worth observing that each $D_\alpha$ is a function if $\Gamma$ 
is complete. In that case each $D_\alpha$ will be countable.

The set $D=\{D_\alpha\mid\alpha\in{\mathcal T}\}$ will be the collection of
domains of the refuting intensional model we are after. Note that, since
each term has a unique type, the $D_\alpha$ are pairwise disjoint. The
$D_\alpha$ are also non-empty. This follows from the assumption that
${\mathcal L}_\alpha\neq\varnothing$ in case $\alpha$ is basic; in case
$\alpha=\langle\alpha_1\ldots\alpha_n\rangle$ it is easy to show that
$\langle\lambda x_{\alpha_1}\ldots\lambda
x_{\alpha_n}.\bot,\varnothing\rangle\in D_\alpha$.

We will define a
function $I$ which will turn out to be an intension function for $D$. First
some handy notation.  If $\pi$ is an ordered pair, write $\pi^1$ and $\pi^2$
for the first and second elements of $\pi$ respectively, so that 
$\pi=\langle\pi^1,\pi^2\rangle$. If $f$ is a function whose values are
ordered pairs, write $f^1$ and $f^2$ for the functions with the same domain
as $f$, such that $f^1(z)=(f(z))^1$ and $f^2(z)=(f(z))^2$ for any argument
$z$. Let $a$ be an assignment for $D$. The substitution $\overleftarrow{a}$
is defined by $\overleftarrow{a}(x)=a^1(x)$ and we let $I^1(a,A)=
A\overleftarrow{a}$ for any term $A$. The second component of $I$, is
defined by letting $I^2=\bigcup_{\alpha\in{\mathcal T}}I^2_\alpha$, where
the $I^2_\alpha$ are functions such that $I^2_\alpha\colon{\mathcal
A}_D\times T^{\mathcal L}_\alpha\to T^{\mathcal L}_\alpha$ if
$\alpha\in{\mathcal B}$ and
\[I^2_\alpha\colon{\mathcal A}_D\times T^{\mathcal L}_\alpha\to
{\mathcal P}(D_{\alpha_1}\times\cdots\times D_{\alpha_n})\ ,\]

\noindent if $\alpha=\langle\alpha_1\ldots\alpha_n\rangle$. The $I^2_\alpha$
in their turn are defined using the following induction on the complexity of
terms.
\begin{enumerate}
\item $I^2_\alpha(a,x_\alpha)=a^2(x)$, if $x$ is a variable;\\
$I^2_\alpha(a,c_\alpha)=c$, if $\alpha$ is basic;\\
$I^2_\alpha(a,c_\alpha)=\{\langle \langle A_1,e_1\rangle,\ldots, 
\langle A_n,e_n\rangle\rangle\mid\langle A_i,e_i\rangle\in D_{\alpha_i}
\mathrel{\&}
\T cA_1\ldots A_n\in\Gamma\}$, 
if $\alpha=\langle\alpha_1\ldots\alpha_n\rangle$;
\item $I^2_{\langle\rangle}(a,\bot)=0$
\item $I^2_{\langle\vec{\alpha}\rangle}
(a,A_{\langle\beta\vec{\alpha}\rangle}B_\beta)= \{\langle\vec{d}\rangle\mid
\langle \langle I^1(a,B),I^2_\beta(a,B)\rangle,
\vec{d}\rangle\in I^2_{\langle\beta\vec{\alpha}\rangle}(a,A)\}$
\item $I^2_{\langle\beta\vec{\alpha}\rangle}
(a,\lambda x_\beta A_{\langle\vec{\alpha}\rangle})=
\{\langle d,\vec{d}\rangle\mid
d\in D_\beta\mathrel{\&}\langle\vec{d}\rangle\in 
I^2_{\langle\vec{\alpha}\rangle}(a[d/x],A)\}$
\item $I^2_{\langle\rangle}(a,A\subset B)=1
\Longleftrightarrow I^2_\alpha(a,A_\alpha)\subseteq I^2_\alpha(a,B_\alpha)$
\end{enumerate}
Note that this definition does not depend on the question whether 
$I$ is an intension function for $D$ and ${\mathcal L}$, and indeed the latter 
is not immediately
obvious. We need to check the requirements in 
definition \ref{intensionfunction}. That $I(a,x)=a(x)$ for any variable
$x$ is immediate and that $I(a,A)=I(a',A)$ if $a$ and $a'$ agree
on the variables free in $A$
follows by a standard property of substitutions and an easy induction. Suppose 
that $B$ is free for $x$ in $A$. Then
\begin{multline*}
        I^1(a,A\{x:=B\})=A\{x:=B\}\overleftarrow{a}=
        A\overleftarrow{a}[x:=B\overleftarrow{a}]=\\
        A\overleftarrow{a}[x:=I^1(a,B)]=
        A\overleftarrow{a[I(a,B)/x]}=I^1(a[I(a,B)/x],A)\ .
\end{multline*}
That $I^2(a,A\{x:=B\})=I^2(a[I(a,B)/x],A)$ if $B$ is free for $x$ in
$A$ follows by a straightforward induction on the complexity of
$A$ which we leave to the reader. 
Thus $I(a,A\{x:=B\})=I(a[I(a,B)/x],A)$ if $B$ is free 
for $x$ in $A$.

It remains to be shown that $I(a,A)\in D_\alpha$
for any assignment $a$ and term $A$
of type $\alpha$. This is done by induction on the complexity of $A$. That
$I(a,x_\alpha)\in D_\alpha$ if $x$ is a variable follows from
the fact that $I(a,x)=a(x)$ and that $I(a,c_\alpha)\in D_\alpha$ if
$\alpha$ is basic is immediate. In the remaining cases the type of $A$ is
complex and it suffices
to prove that whenever $\alpha=\langle\alpha_1\ldots\alpha_n\rangle$,
$\langle B_1,e_1\rangle\in D_{\alpha_1},\ldots,$ and $
\langle B_n,e_n\rangle\in D_{\alpha_n}$:
\begin{itemize}
\item[(a)] If $\T A\overleftarrow{a}B_1\ldots B_n\in\Gamma$ then 
$\langle \langle B_1,e_1\rangle,\ldots,\langle B_n,e_n\rangle\rangle
\in I^2(a,A)$;
\item[(b)] If $\F A\overleftarrow{a}B_1\ldots B_n\in\Gamma$ then 
$\langle \langle B_1,e_1\rangle,\ldots,\langle B_n,e_n\rangle\rangle
\notin I^2(a,A)$.
\end{itemize}
We shall consider each case. IH will be short for `induction 
hypothesis'.
\begin{itemize}
\item $A_\alpha\equiv c$ and $\alpha=\langle\alpha_1\ldots\alpha_n\rangle$. 
The requirement follows from the definition of $I^2(a,c)$ and clause 1.\ of
definition \ref{hintikka}.
\item $A_\alpha\equiv \bot$ and $\alpha=\langle\rangle$. The (a) part of the
property follows from clause 2.\ of definition \ref{hintikka}, the (b) part
from the fact that $I^2(a,\bot)=0=\varnothing$.
\item $A\equiv B_{\langle\beta\alpha_1\ldots\alpha_n\rangle}C_\beta$. Suppose 
$\langle B_1,e_1\rangle\in D_{\alpha_1},\ldots,
\langle B_n,e_n\rangle\in D_{\alpha_n}$, then
\begin{eqnarray*}
\T (BC)\overleftarrow{a}B_1\ldots B_n\in\Gamma&\Longleftrightarrow&\\
\T B\overleftarrow{a}C\overleftarrow{a}B_1\ldots B_n\in\Gamma&\Longrightarrow&
\mbox{(IH)}\\
\langle I(a,C),\langle B_1,e_1\rangle,\ldots,\langle B_n,e_n\rangle\rangle
\in I^2(a,B)&\Longleftrightarrow&\mbox{(def.\ of $I$)}\\
\langle\langle B_1,e_1\rangle,\ldots,\langle B_n,e_n\rangle\rangle
\in I^2(a,BC)
\end{eqnarray*}
This proves the (a) part of the case; the (b) part is similar.
\item $A\equiv (\lambda x_{\alpha_1}C_{\langle\alpha_2\ldots\alpha_n\rangle})$. 
Again suppose $d_1=\langle B_1,e_1\rangle\in D_{\alpha_1},\ldots,
d_n=\langle B_n,e_n\rangle\in D_{\alpha_n}$, and reason as follows.
\begin{eqnarray*}
\F (\lambda x.C)\overleftarrow{a}B_1\ldots B_n\in\Gamma&\Longleftrightarrow&\\
\F \lambda x.C\overleftarrow{a}_xB_1\ldots B_n\in\Gamma&\Longrightarrow&
\mbox{Def.\ \ref{hintikka}, $B_1$ is closed}\\
\F C\overleftarrow{a}_x\{x:=B_1\}B_2\ldots B_n\in\Gamma&\Longleftrightarrow&\\
\F C\overleftarrow{a[d_1/x]}B_2\ldots B_n\in\Gamma&\Longrightarrow&\mbox{(IH)}\\
\langle d_2,\ldots,d_n\rangle
\notin I^2(a[d_1/x],C)&\Longleftrightarrow&\mbox{(def.\ of $I$)}\\
\langle d_1,d_2,\ldots,d_n\rangle\notin I^2(a,\lambda x.C)
\end{eqnarray*}
This proves the (b) part, which is similar to the (a) part.
\item $A_\alpha\equiv B\subset C$. Then $\alpha=\langle\rangle$ and
$B$ and $C$ have some type $\langle\alpha_1\ldots\alpha_n\rangle$. 
Using induction we
may assume that $I(a,B),I(a,C)\in D_{\langle\alpha_1\ldots\alpha_n\rangle}$.
Suppose $\T (B\subset C)\overleftarrow{a}\in\Gamma$, i.e.\
$\T B\overleftarrow{a}\subset C\overleftarrow{a}\in\Gamma$ and
reason as follows.
\begin{eqnarray*}
\langle \langle B_1,e_1\rangle,\ldots,
\langle B_n,e_n\rangle\rangle\in I^2(a,B)&\Longrightarrow& \mbox{(IH)}\\
\F B\overleftarrow{a}B_1\ldots B_n\notin\Gamma&\Longrightarrow& 
\mbox{(Def.\ \ref{hintikka})}\\
\T C\overleftarrow{a}B_1\ldots B_n\in\Gamma&\Longrightarrow& \mbox{(IH)}\\
\langle \langle B_1,e_1\rangle,\ldots,
\langle B_n,e_n\rangle\rangle\in I^2(a,C)
\end{eqnarray*}
We conclude that $I^2(a,B)\subseteq I^2(a,C)$ and that 
$I^2(a,B\subset C)=1$. This
proves the (a) part of the property. The (b) part is left to
the reader.
\end{itemize}
This concludes the
proof that $I$ is an intension function for $D$ and ${\mathcal L}$. 
Now define the function $E$ by
letting $E(\langle A,e\rangle)=e$ if 
$\langle A,e\rangle\in D_\alpha$ for any
complex $\alpha$. Clearly, $E(I(a,A))=I^2(a,A)$ for any $A_\alpha$, 
$E$ is an extension function for $D$,
and $M=\langle D,I,E\rangle$ is an intensional model for the 
language ${\mathcal L}$.
It is easy to see that $M$ refutes $\Gamma$. We have already established
that $M$ is countable if $\Gamma$ is complete, and proposition \ref{normalize}
gives a normal intensional model refuting $\Gamma$ which is countable in case
$\Gamma$ is complete.

\end{proof}
\noindent Before we continue with the proof of Model Existence, let us look
at an application. Hintikka's Lemma sometimes gives an easy way of showing the
refutability of certain sequents. For example, while standard higher order
logic validates the sentence $p=q\lor q=r\lor r=p$ (where $p$, $q$ and $r$
are type $\langle\rangle$ constants), a Hintikka sequent (corresponding to an
open tableau branch) shows that this is not the case in \textsf{ITL}.
\begin{proposition}
$\Rightarrow p=q,q=r,r=p$ is refutable by an intensional model.
\end{proposition}
\begin{proof}
The non-abbreviated form of $\Rightarrow p=q,q=r,r=p$ is\\[1ex]
$\Gamma=\{\F(\lambda z.\top)\subset(\lambda z.zp\subset zq),\
\F(\lambda z.\top)\subset(\lambda z.zq\subset zr),$\\ \hspace*{\fill}$
\F(\lambda z.\top)\subset(\lambda z.zr\subset zp)\}$,\\[1ex]
with $z$ a
variable of type $\langle\langle\rangle\rangle$.
An inspection tells that the following extension $\Gamma^+$ of $\Gamma$, in which 
$c_1$, $c_2$ and $c_3$
are constants of type $\langle\langle\rangle\rangle$, is a Hintikka sequent.\\[1ex]
$\{\F(\lambda z.\top)\subset(\lambda z.zp\subset zq),\ 
\T(\lambda z.\top)c_1,\ 
\F(\lambda z.zp\subset zq)c_1,\ 
\T\top,$\\$ 
\F c_1p\subset c_1q,\ 
\T c_1p,\ 
\F c_1q,\
\F(\lambda z.\top)\subset(\lambda z.zq\subset zr),\
\T(\lambda z.\top)c_2,$\\$
\F(\lambda z.zq\subset zr)c_2,\
\F c_2q\subset c_2r,\
\T c_2q,\
\F c_2r,$\\$
\F(\lambda z.\top)\subset(\lambda z.zr\subset zp),\
\T(\lambda z.\top)c_3,\
\F(\lambda z.zr\subset zp)c_3,\
\F c_3r\subset c_3p,$\\$
\T c_3r,\
\F c_3p\}
$\\[1ex]
It follows that $\Gamma^+$ and hence $\Gamma$ are refutable by an
intensional model. 
\end{proof}

\noindent This shows that it is consistent to assume that there are at least
three propositions. It is clear that the method can be generalized to
show that it is consistent to assume a set of propositions $\geq$ any given
countable cardinality. If the countability restriction on languages is dropped,
the existence of intensional models with type $\langle\rangle$ domains
$\geq$ any given cardinality is obtained. We leave it to the reader to show
that $p\leftrightarrow q\Rightarrow p=q$ and other instances of
Extensionality are refutable.

Let us return to the main line of argument. In order to state the model
existence theorem below, we need the notion of a \emph{provability property}
(closely related to Smullyan's~\cite{smul:firs68} \emph{abstract consistency
property}).

\begin{definition}
Let ${\mathcal P}$ be a set of sequents in the language ${\mathcal L}$.
${\mathcal P}$ is a \emph{provability property} in ${\mathcal L}$
if ${\mathcal P}$ is closed under sequent rules, i.e.\ if
$\Gamma\in{\mathcal P}$
whenever $\{\Gamma_1,\ldots,\Gamma_n\}\subseteq{\mathcal P}$ and
$\Gamma_1,\ldots,\Gamma_n/\Gamma$ is a sequent rule. 

A provability
property ${\mathcal P}$ in ${\mathcal L}$ is \emph{sound} if no 
$\Gamma\in{\mathcal P}$
is refuted by an intensional model for ${\mathcal L}$.
\end{definition}

\noindent We now come to Model Existence itself: sequents that are
\emph{not} elements of a sound provability property (in an extended language)
can be extended to Hintikka sequents (in that language) and are hence
refutable.
\begin{theorem}[Model Existence]\label{modelexistence}
Let ${\mathcal L}$ and ${\mathcal C}$ be languages 
such that ${\mathcal L}\cap{\mathcal C} =\varnothing$ and
each ${\mathcal C}_\alpha$ is denumerably infinite. Assume that ${\mathcal
P}$ is a sound provability property in
${\mathcal L}\cup{\mathcal C}$ and that  $\Gamma$ is a sequent in the
language ${\mathcal L}$. If $\Gamma \notin {\mathcal P}$ then $\Gamma$ is
refuted by a countable normal intensional model.
\end{theorem}%
\begin{proof}
Let ${\mathcal P}$ and $\Gamma$ be as described. We construct a
Hintikka sequent $\Gamma^*$ such that $\Gamma \subseteq
\Gamma^*$. Let $\vartheta_1,\ldots,\vartheta_n,\ldots$ be an enumeration of 
all signed sentences in ${\mathcal L}\cup{\mathcal C}$. Write $\iota(\vartheta)$
for the index that the signed
sentence $\vartheta$ obtains in this enumeration. Let $\Gamma_0=\Gamma$ and
define each $\Gamma_{n+1}$ by distinguishing the following cases.
\begin{itemize}
\item $\Gamma_{n+1}=\Gamma_n$,
        if $\Gamma_n \cup \{\vartheta_n\} \in {\mathcal P}$;
\item $\Gamma_{n+1}=\Gamma_n \cup \{\vartheta_n\}$,
        if $\Gamma_n \cup \{\vartheta_n\} \notin {\mathcal P}$ and $\vartheta_n$ is
        not of the form $\F A\subset B$;
\item $\Gamma_{n+1}=\Gamma_n \cup 
\{\vartheta_n, \T Ac_1\ldots c_n,\F Bc_1\ldots c_n\}$,
        if $\Gamma_n \cup \{\vartheta_n\} \notin {\mathcal P}$ and
$\vartheta_n=\F A\subset B$ for $A$ and $B$ of
type $\langle\alpha_1\ldots\alpha_n\rangle$, 
where each $c_i$ is the first constant in ${\mathcal C}_{\alpha_i}$
which does not occur in $\Gamma_n \cup \{\vartheta_n\}$ and is no
element of $\{c_1,\ldots,c_{i-1}\}$
\end{itemize}
This is well-defined since each $\Gamma_n$ contains only a finite
number of constants from ${\mathcal C}$. That $\Gamma_n\notin{\mathcal P}$ for each
$n$ follows by a simple induction which uses the definition of a
provability property and the fact that $[\mathop{\subset}{\sf R}]$ is a
sequent rule. Define $\Gamma^* = \bigcup_n \Gamma_n$.
For all finite sets
$\{\vartheta_{k_1},\ldots,\vartheta_{k_n}\}$ and for all $k \ge max\{k_1,\ldots,k_n\}$
\begin{equation}\label{fin}
\{\vartheta_{k_1},\ldots,\vartheta_{k_n}\}\subseteq\Gamma^* \Leftrightarrow
\Gamma_k \cup \{\vartheta_{k_1},\ldots,\vartheta_{k_n}\} \notin {\mathcal P}
\end{equation}
In order to show this, let $k \ge max\{k_1,\ldots,k_n\}$
and let
$\{\vartheta_{k_1},\ldots,\vartheta_{k_n}\}\subseteq\Gamma^*$. Then there is
some $\ell$ such that
$\{\vartheta_{k_1},\ldots,\vartheta_{k_n}\}\subseteq\Gamma_\ell$. Let
$m=max\{k,\ell\}$. We have that $\Gamma_k \cup
\{\vartheta_{k_1},\ldots,\vartheta_{k_n}\}\subseteq\Gamma_m$. Since
$\Gamma_m\notin{\mathcal P}$ and ${\mathcal P}$ is closed under supersets (rule $[W]$),
it
follows that $\Gamma_k \cup \{\vartheta_{k_1},\ldots,\vartheta_{k_n}\}
\notin {\mathcal P}$. For the reverse direction, suppose that $\Gamma_k
\cup \{\vartheta_{k_1},\ldots,\vartheta_{k_n}\} \notin {\mathcal P}$. Then,
since ${\mathcal P}$ is closed under supersets, $\Gamma_{k_i}
\cup \{\vartheta_{k_i}\} \notin {\mathcal P}$, for each of the $k_i$. By the
construction of $\Gamma^*$ each $\vartheta_{k_i}\in\Gamma^*$ and
$\{\vartheta_{k_1},\ldots,\vartheta_{k_n}\}\subseteq\Gamma^*$. 

With the help of (\ref{fin}) it can be verified that $\Gamma^*$ is a
Hintikka sequent. The last condition of Definition
\ref{hintikka} immediately follows from the construction of 
$\Gamma^*$. We check condition 5, which may serve as an example
for the other cases. Assume $\T A\subset B\in \Gamma^*$ and
let $k$ be the maximum of $\iota(\T A\subset B)$, 
$\iota(\T B\vec{C})$, and $\iota(\F A\vec{C})$. Since, by (\ref{fin}),
$\Gamma_k\cup\{\T A\subset B\}\notin{\mathcal P}$ and
since ${\mathcal P}$ is closed under sequent rules, it must be the case
that either $\Gamma_k\cup\{\T B\vec{C}\}\notin{\mathcal P}$ or
$\Gamma_k\cup\{\F A\vec{C}\}\notin{\mathcal P}$, Using (\ref{fin}),
we find that $\T B\vec{C}\in \Gamma^*$ or $\F A\vec{C}\in \Gamma^*$.

We conclude that $\Gamma^*$ is refuted by an intensional model $M$. In
order to prove that there is a normal countable intensional model
that refutes $\Gamma^*$ and hence $\Gamma$ it suffices to show
that $\Gamma^*$ is complete. Let $\varphi$ be any sentence of
${\mathcal L}\cup{\mathcal C}$ and assume that $\T\varphi\notin\Gamma^*$
and $\F\varphi\notin\Gamma^*$. Then, by (\ref{fin}),
$\Gamma_k\cup\{\T\varphi\}\in{\mathcal P}$ and 
$\Gamma_k\cup\{\F\varphi\}\in{\mathcal P}$, for sufficiently large $k$.
But $M$ refutes $\Gamma_k$ and therefore must either refute
$\Gamma_k\cup\{\T\varphi\}$ or $\Gamma_k\cup\{\F\varphi\}$, contradicting
the soundness of ${\mathcal P}$. Thus $\Gamma^*$ is complete and
some normal countable intensional model refutes $\Gamma^*$ and $\Gamma$.

\end{proof}
\noindent From model existence we can derive some nice corollaries. In the
following $\Gamma$ will always be a sequent in some language ${\mathcal L}$
while $\Delta$ ranges over sequents in ${\mathcal L}\cup {\mathcal C}$,
where ${\mathcal L}$ and ${\mathcal C}$ are as in the formulation of Theorem
\ref{modelexistence}.
\begin{corollary}[Generalized Compactness]
If $\Gamma$ is i-valid
then some finite $\Gamma_0 \subseteq \Gamma$ is
i-valid.
\end{corollary}
\begin{proof}
$\{\Delta\mid$ some finite $\Delta_0 \subseteq \Delta$ is
i-valid$\}$ is a sound provability property. 
\end{proof}
\begin{corollary}[Generalized L\"owenheim--Skolem]
If $\Gamma$ is not i-valid then $\Gamma$ is refutable by a countable normal
intensional model.
\end{corollary}
\begin{proof}
$\{\Delta\mid\Delta$ is i-valid$\}$ is a sound provability property.
\end{proof}
\begin{corollary}[Generalized Completeness]
If $\Gamma$ is i-valid then $\Gamma$ is provable. Hence 
$\Pi\models_i\Sigma\Longrightarrow\Pi\vdash\Sigma$.
\end{corollary}
\begin{proof}
$\{\Delta\mid\Delta$ is provable$\}$ is a sound provability
property.
\end{proof}

\begin{corollary}[Cut elimination]
If $\Pi,\varphi\vdash\Sigma$ and $\Pi\vdash\Sigma,\varphi$ then
$\Pi\vdash\Sigma$.
\end{corollary}
\begin{proof}
Use soundness and completeness.
\end{proof}

\section{A Linguistic Application}\label{ling}
We now turn to a linguistic application of \textsf{ITL} and will develop the
semantics of a tiny fragment of English containing verbs of propositional
attitude. It will be shown that, given the present logic, it is consistent
for an agent $a$ to know that $\varphi$ without knowing that $\psi$, even if
$\varphi$ and $\psi$ are co-entailing.

Before considering our special application, however, let us address the
general point of axiomatic extensions of the base logic. In most
applications one will like to work with a subclass of the class of intensional
models that conform to some set of non-logical axioms ${\mathcal S}$. In
that case one can define $\Pi\models_{\mathcal S}\Sigma$ to be ${\mathcal
S}\cup\Pi\models_i\Sigma$, while $\Pi\vdash_{\mathcal S}\Sigma$ can be
defined as ${\mathcal S}\cup\Pi\vdash\Sigma$. Soundness and generalized
completeness immediately give that $\Pi\models_{\mathcal
S}\Sigma\Longleftrightarrow\Pi\vdash_{\mathcal S}\Sigma$. Not all
applications will instantiate ${\mathcal S}$ in the same way, but one set of
axioms that immediately come to mind, and that we shall adopt here, are the
usual principles of $\lambda$-conversion. We may add these by assuming that
${\mathcal S}$ contains all universal closures of instantiations of the
following schemes.
\begin{enumerate}
  \item[($\alpha$)] $\lambda x.A=\lambda y.A\{x:=y\}$, if $y$ is free for 
  $x$ in $A$;
  \item[($\beta$)] $(\lambda x.A)B=A\{x:=B\}$, if $B$ is free for $x$ in $A$;
  \item[($\eta$)] $\lambda x.Ax=A$, if $x$ is not free in $A$.
\end{enumerate}
As soon as these schemes are added to the base logic, the result is full
intensional identity of $\beta\eta$ equivalent terms, i.e.\ 
$\models_{\mathcal S} A=B$ will hold if $A=_{\beta\eta}B$.

\begin{table}[t]
\begin{tabular*}{\textwidth}%
{@{\qquad}@{\extracolsep{\fill}}rl@{\qquad\qquad}rl}\hline
{\sc word}&{\sc translation}&{\sc word}&{\sc translation}\\ \hline
if        &$\lambda p_{\langle\rangle}\lambda q_{\langle\rangle}.p\to q$&
man       &$\mbox{\it man}_{\langle e\rangle}$\\
no        &$\lambda P'_{\langle e\rangle}\lambda P_{\langle e\rangle}.
           \neg\exists x_e(P'x\land Px)$&
unicorn   &$\mbox{\it unicorn}_{\langle e\rangle}$\\
some      &$\lambda P'_{\langle e\rangle}\lambda P_{\langle e\rangle}.
           \exists x_e(P'x\land Px)$&
runs     &$\mbox{\it run}_{\langle e\rangle}$\\	   
every     &$\lambda P'_{\langle e\rangle}\lambda P_{\langle e\rangle}.
           \forall x_e(P'x\to Px)$&
laughs       &$\mbox{\it laugh}_{\langle e\rangle}$\\	   
loves      &$\lambda Q_{\langle\langle e\rangle\rangle}\lambda x_e.
           Q(\lambda y_e.\mbox{\it love}_{\langle ee\rangle}\,xy)$&
Bill       & $\mbox{\it bill}_{\langle\langle e\rangle\rangle}$\\   
is    &$\lambda Q_{\langle\langle e\rangle\rangle}\lambda x_e.
           Q(\lambda y_e.x=y)$&
Ann        & $\mbox{\it ann}_{\langle\langle e\rangle\rangle}$\\	   
knows     &$\lambda p_{\langle\rangle}\lambda x_e.
           \mbox{\it know}_{\langle e\langle\rangle\rangle}\,xp$&
Tully       & $\mbox{\it tully}_{\langle\langle e\rangle\rangle}$\\	   
believes  &$\lambda p_{\langle\rangle}\lambda x_e.
           \mbox{\it believe}_{\langle e\langle\rangle\rangle}\,xp$&
Cicero       & $\mbox{\it cicero}_{\langle\langle e\rangle\rangle}$\\	   
\hline
\end{tabular*}
\caption{\label{lex}Some words and their translations}
\end{table}

For our linguistic application we will proceed along lines pioneered by
Montague~\cite{mont:prop73} and define a small fragment of English. The
words of this fragment are given in Table \ref{lex}, along with their
translations into type logic. In these translations the terms \textit{love},
\textit{run}, \textit{man}, etc.\ are constants of the types indicated, where
$e$ is the type of \emph{entities}. The
set of \emph{syntactic structures} is obtained by stipulating that all words
in Table \ref{lex} are syntactic structures and that $[XY]$ is a syntactic
structure whenever $X$ and $Y$ are syntactic structures. Defining syntactic
structures in this way leads to a lot of gibberish along with the
structures we are interested in, but this is not important for present
purposes. As long as the desired structures are there and get reasonable
interpretations our aim is served.

Let us define the relation $\leadsto$ (``translates as'') between syntactic
structures and terms as the smallest relation such that 1) $X\leadsto A$ if
$X$ is a word and $A$ is its translation in Table \ref{lex} and 2) if
$X\leadsto A$ and $Y\leadsto B$ then $[XY]\leadsto AB$ if $AB$ is a
well-formed term and $[XY]\leadsto BA$ if $BA$ is well-formed. This leaves
open the possibility that a syntactic structure does not get a translation
and indeed many do not. Structures $X$ for which there is no $A$ such
that $X\leadsto A$ are called \emph{uninterpretable} and we have no interest
in them.

Let us turn to some syntactic structures that are interpretable. In
(\ref{noman}) below two are given, together with (the $\beta$ normal forms
of) their interpretations. Clearly, (\ref{noman}b), the interpretation of
(\ref{noman}a), i-entails and is i-entailed by (\ref{noman}d), which is the
interpretation of (\ref{noman}c).

\eenumsentence{\label{noman}
  \item {}[[[no man]laughs][if[[some unicorn]runs]]]
  \item $\exists x(\mbox{\it unicorn}\;x\land\mbox{\it run}\;x)\to
        \lnot\exists x(\mbox{\it man}\;x\land\mbox{\it laugh}\;x)$
  \item {}[[[no unicorn]runs][if[[some man]laughs]]]
  \item $\exists x(\mbox{\it man}\;x\land\mbox{\it laugh}\;x)\to
        \lnot\exists x(\mbox{\it unicorn}\;x\land\mbox{\it run}\;x)$
}
This does not mean however that (\ref{noman}b) and (\ref{noman}d) are
identical in all intensional models, as nothing excludes the possibility 
that $I(a,(\ref{noman}b))\neq I(a,(\ref{noman}d))$ for some intension
function $I$. It follows that the two structures in (\ref{knowsnoman}) are
\emph{not} co-entailing.
\eenumsentence{\label{knowsnoman}
  \item {}[[every man][knows[[[no man]laughs][if[[some unicorn]runs]]]]]
  \item $\forall y(\mbox{\it man}\;y\to\mbox{\it know}\;y\,
         (\exists x(\mbox{\it unicorn}\;x\land\mbox{\it run}\;x)\to
        \lnot\exists x(\mbox{\it man}\;x\land\mbox{\it laugh}\;x)))$
  \item {}[[every man][knows[[[no unicorn]runs][if[[some man]laughs]]]]]
  \item $\forall y(\mbox{\it man}\;y\to\mbox{\it know}\;y\,
        (\exists x(\mbox{\it man}\;x\land\mbox{\it laugh}\;x)\to
        \lnot\exists x(\mbox{\it unicorn}\;x\land\mbox{\it run}\;x)))$
}
Suppose that $c$ is some constant of type $e$. Then 
$\mbox{\it know}\;c\,(\ref{noman}b)\Rightarrow\mbox{\it
know}\;c\,(\ref{noman}d)$ is in fact a Hintikka sequent and is therefore
refuted by an intensional model (addition of ($\alpha$), ($\beta$) and ($\eta$)
does not change this). This intensional model can also be used to show that
(\ref{knowsnoman}b) does not entail (\ref{knowsnoman}d). This is as desired,
for even if (\ref{knowsnoman}a) holds there may well be a man who has not
managed to draw the inference necessary to arrive at (\ref{noman}c).
We have thus shown that the logic avoids the problem of logical omniscience
in the sense that it does not exclude the possibility that a person knows
one thing but fails to know another thing logically equivalent with it.
Essential use was made of the failure of Extensionality in our logic
\textsf{ITL}: terms of complex type can have the same extensions, even in
all intensional models, without necessarily having the same intension.

This distinction between extension and intension does not extend to terms of
basic type however and this raises the question how \emph{names} are to be
dealt with. If they are treated straightforwardly using constants of type
$e$ (e.g.\ $b_e$, or in the present context preferably $\lambda P.Pb$, for
`Bill') we run into the standard problems of the `Cicero--Tully' or
`Hesperus--Phosphorus' kind. However, there are many reasonable translations
that do not directly equate names with type $e$ constants. The translations
in Table \ref{lex}, that send names to constants of the quantifier type
$\langle\langle e\rangle\rangle$, may serve as an example, provided some
meaning postulates (additions to ${\mathcal S}$) like the following are
adopted.

\eenumsentence{\label{mpnames}
\item $\forall P(\mbox{\it ann}\,P\leftrightarrow Pa)$
\item $\forall P(\mbox{\it bill}\,P\leftrightarrow Pb)$
\item $\forall P(\mbox{\it tully}\,P\leftrightarrow Pt)$
\item $\forall P(\mbox{\it cicero}\,P\leftrightarrow Pc)$
}
The structure [Tully runs] translates as $\mbox{\it
tully}\;\mbox{\it run}$, but given the meaning postulates just introduced,
this is equivalent with $\mbox{\it run}\;t$. Similarly, [Cicero runs]
translates as $\mbox{\it cicero}\;\mbox{\it run}$, equivalent with
$\mbox{\it run}\;c$. And since [Tully[is Cicero]]
is translated as $\mbox{\it tully}(\lambda x.\mbox{\it cicero}(\lambda
y.x=y))$, which is equivalent with $t=c$, it is readily explained why the
argument \emph{Tully runs, Tully is Cicero, therefore Cicero runs} holds.
But this reasoning essentially depends on extensional equivalence and
therefore will not go through once propositional attitudes enter the picture.
Consider the structure [Ann[believes[Tully runs]]]. It translates as 
$\mbox{\it ann}(\lambda x.\mbox{\it believe}\;x(\mbox{\it
tully}\;\mbox{\it run}))$ and this is equivalent with
$\mbox{\it believe}\;a\,(\mbox{\it tully}\;\mbox{\it run})$, while 
$\mbox{\it believe}\;a\,(\mbox{\it cicero}\;\mbox{\it run})$
is equivalent with the
translation of [Ann[believes [Cicero runs]]]. However, there is no
co-entailment between these sentences, even in the presence of the
postulates in (\ref{mpnames}) and the translation of [Tully[is Cicero]].

This shows that even for names the sense/reference distinction can be
captured in this logic, provided one is willing to treat names with the help
of predicates (Quine's `primacy of predicates' comes to mind). Treating
them as being of type $\langle\langle e\rangle\rangle$, as we have done
here, is one possible strategy. There may be others.

The present application of our intensional type theory to linguistic
semantics has avoided the concept of possible worlds altogether, as it was
not needed in order to illustrate our points. However, as possible worlds
are obviously extremely useful for the analysis of a range of natural
language constructions (though not for the true intensionality we have been
concerned with in this paper), one might well want to combine them with the
present approach. Muskens~\cite[chapter 4]{mp} gives a translation of what
is essentially the fragment of Montague~\cite{mont:prop73} into a two-sorted
relational type theory, with possible worlds providing an additional basic
type. Although the type theory in~\cite{mp} validates Extensionality, its
language essentially is the language employed here, so that the translation
can also serve as a translation into \textsf{ITL}. A minor variation
will treat names as they are treated above.

\section{Worlds}\label{worlds}
\textsf{ITL} is a generalization of the usual formulation of type theory and
intensionality is obtained by giving up the axiom of Extensionality, not by
the introduction of possible worlds, as in modal logic. However, while the
usual Kripke-style semantics is not known to do a very good job regarding
the puzzles of intensionality we have been concerned with here, it does
perform very well when it comes to modal reasoning, temporal reasoning,
counterfactual reasoning etc. So it seems that worlds and the possibility to
quantify over worlds are still welcome, even to those who accept the claim
that the present approach to intensionality is superior to the modal one.

If such a combination of modality with true intensionality is desired, one
way to proceed would be to simply add domains of worlds to the existing
intensional models and interpret a modal higher order language on the
results, a course of action followed in Muskens~\cite{homl}. There is,
however, an easier way. Once true intensionality is obtained in the way it
was done in this paper, worlds can also be constructed out of propositions,
the inhabitants of the domain $D_{\langle\rangle}$, while accessibility
relations can be obtained as well. The procedure will only be sketched in this
section; more formal considerations and comparisons with standard
approaches to modality will be left to a future occasion.

The idea of constructing possible worlds out of other entities is an old
one. E.g.\ Wittgenstein~\cite{witt:trac22} constructs them out of `states of
affairs' and Carnap~\cite{carn:mean} takes worlds to be
`state-descriptions', maximal consistent sets of sentences. A recent
construction of worlds from propositions can be found in
Pollard~\cite{poll:hype05}. Varying upon such proposals, one can identify
worlds with certain objects of type $\langle\langle\rangle\rangle$ here,
i.e.\ objects whose extensions are sets of propositions. Here is the
construction. Assume, for simplicity, the principles $(\alpha)$, $(\beta)$,
and $(\eta)$ discussed above, and let $\Omega$, which will stand for the
predicate `is a world', be a fixed constant of type
$\langle\langle\langle\rangle\rangle\rangle$, while $w$ varies over objects
of type $\langle\langle\rangle\rangle$. Stipulate the following.

\begin{enumerate}
\item[(W1)] $\forall w(\Omega w\to\lnot w\bot)$
\item[(W2)] $\forall w(\Omega w\to (w(A\subset B)\leftrightarrow
             \forall\vec{x}(w(A\vec{x})\to w(B\vec{x}))))$
\end{enumerate}
The first of these axioms requires world extensions to be consistent while
addition of the second schema makes worlds distribute over logical
operators. Statements such as the following become derivable.

\begin{enumerate}
\item[a.] $\forall w(\Omega w\to (w(\lnot\varphi)
          \leftrightarrow\lnot(w\varphi)))$
\item[b.] $\forall w(\Omega w\to (w(\varphi\land\psi)
          \leftrightarrow((w\varphi)\land(w\psi))))$
\item[c.] $\forall w(\Omega w\to (w(\forall
          x\varphi)\leftrightarrow\forall x(w\varphi)))$
\item[d.] $\forall w(\Omega w\to (w(\exists
          x\varphi)\leftrightarrow\exists x(w\varphi)))$
\end{enumerate}
The first of these statements says that worlds are complete, while the last
two are `Henkin properties' that enforce, for example, that if an
existential proposition is an element of the extension of a given world some
proposition witnessing the existential must also be an element. In general,
given (W1) and (W2), worlds single out sets of propositions that could be
simultaneously true. 

It is natural from this perspective to introduce a constant $w_0$ of type
$\langle\langle\rangle\rangle$ that is meant to denote the \emph{actual
world}, the world consisting of all propositions that \emph{are} true (in a
given intensional model). If this is wanted one should stipulate the
following.

\begin{enumerate}
\item[a.] $\Omega w_0$
\item[b.] $\forall p_{\langle\rangle}(w_0p\leftrightarrow p)$. 
\end{enumerate}
The first of these requirements merely stipulates that $w_0$ is a world
while the second makes it the actual world. In models additionally
satisfying Extensionality, $w_0$ is the \emph{only} world in view of the
fact that there are exactly two propositions in such models, but in other
intensional models there is no such trivialization. Note, by the way, that
nothing requires co-extensional worlds to be identical. The set of
propositions that are true in some world does not necessarily determine it.

Since worlds are of type $\langle\langle\rangle\rangle$ it is possible to
iterate and form propositions $w_1\varphi$, $w_2(w_1\varphi)$,
$w_3(w_2(w_1\varphi))$ and so on. Is it acceptable that these differ in
truth value? Here we shall assume that this is not so and that if a
proposition $w\varphi$ is true at some world, it is true at all,
provided $w$ is a world. The question whether a proposition is in the
extension of a world should arguably be world-independent itself. The way to
enforce this is by adopting the following axiom scheme.

\begin{enumerate}
\item[(W3)] $\forall ww'((\Omega w\land\Omega w')\to (w(w'\varphi)
            \leftrightarrow (w'\varphi)))$
\end{enumerate}
In a similar vein, the question whether an $\langle\langle\rangle\rangle$
object is a world, presumably should also be uniform across worlds:
\begin{enumerate}
\item[(W4)] $\forall w(\Omega w\to\forall w'(\Omega w'\leftrightarrow
w(\Omega w')))$
\end{enumerate}
We now have worlds, but we still do not have accessibility relations between
worlds. These can be obtained, however, by considering more expressions of
type $\langle\langle\langle\rangle\rangle\rangle$. If $R$ is such an
expression, it can be interpreted as the predicate `is accessible', and
$\lambda w\lambda w'.w(Rw')$ will play the role of an accessibility
relation. The usual relational properties (transitivity, reflexivity,
euclideanness,\ldots) can then either be stipulated or, depending on the
choice of $R$, be shown to hold. For example, the universal accessibility
relation $\lambda w\lambda w'.w(\Omega w')$ is easily seen to be an
equivalence relation on the set of worlds in view of (W4).

A next step is the introduction of the usual modal operators. Modal boxes
can be obtained by writing $[R]$ for $\lambda p\forall w((\Omega w\land
Rw)\to wp)$, so that $[R]\varphi$ will reduce to $\forall w((\Omega w\land
Rw)\to w\varphi)$. Diamonds are obtained as usual, as the duals of boxes:
$\langle R\rangle$ is short for $\lambda p.\lnot[R]\lnot p$. Note that if
$w'$ can be shown to be a world, the statement $w'([R]\varphi)$, i.e.\
$w'(\forall w((\Omega w\land Rw)\to w\varphi))$, will be equivalent with
$\forall w((w'(\Omega w)\land w'(Rw))\to w'(w\varphi))$ by the distribution
of worlds over logical operators and the last statement will in its turn be
equivalent with $\forall w((\Omega w\land w'(Rw))\to w\varphi)$ by (W3) and
(W4).

Let us give another example of an accessibility relation some of whose
properties follow from its definition. The relation of belief considered in
the previous section, $\lambda p\lambda x.\mbox{\it believe}\,xp$, is one of
\emph{explicit} belief. It is not closed under entailment or even under
logical equivalence. But there is also a notion of \emph{implicit} belief
that is closed under entailment. Roughly, one implicitly believes $\varphi$
if one rationally \emph{should} believe $\varphi$ given one's explicit
beliefs. One way to model this (for some arbitrary agent \emph{john}) is to
consider the following property $R$ of worlds. 
\[\lambda w.\forall p((\mbox{\it believe}\,\mbox{\it john}\,p
\leftrightarrow w(\mbox{\it believe}\,\mbox{\it john}\,p))\land
(\mbox{\it believe}\,\mbox{\it john}\,p\to wp))\]
Here a world $w$ is accessible if John's explicit beliefs in $w$ are exactly
those that John actually holds and if those explicit beliefs are in fact
true in $w$. There may fail to be such worlds, for example if John's
explicit beliefs are in fact inconsistent, a situation not ruled out by our
previous considerations. But it is possible to stipulate that 
$\lambda w\lambda w'.w(Rw')$ is in fact serial:
\[\forall w\exists w'\forall p((w(\mbox{\it believe}\,\mbox{\it john}\,p)
\leftrightarrow w'(\mbox{\it believe}\,\mbox{\it john}\,p))\land
(w(\mbox{\it believe}\,\mbox{\it john}\,p)\to w'p))\]
Such a stipulation is simultaneously an existence requirement on worlds and
a rationality constraint on John's beliefs. It will lead to the derivability
of the usual \textbf{D} axiom, as $[R]\varphi\to\langle R\rangle\varphi$
will now hold for all $\varphi$. Note that the definition of $R$ immediately
gives transitivity and euclideanity of $\lambda w\lambda w'.w(Rw')$, so that
we have, for all $\varphi$, that $[R]\varphi\to[R][R]\varphi$ and $\langle
R\rangle\varphi\to[R]\langle R\rangle\varphi$. These correspond to the usual
\textbf{4} and \textbf{5} axioms.

\section{Conclusion}
In this paper we have introduced an abstract and simple notion of
\emph{intensional model} that is a generalization of Henkin's general
models. Its definition does not involve concepts that have no immediate
intuitive justification, such as the ``abstraction designation functions''
of Fitting~\cite{fitt:type02} or the ``application operators'' of 
Benzm\"uller et al.~\cite{BBK04}. These operators provide generalized,
non-standard notions of abstraction in one case and of application in the
other, but seem to have no justification other than a purely technical one.
The present approach, in contrast, gives a kind of minimal logic of
intension and extension, with ingredients that well-nigh \emph{any} logic of
intension and extension seems to need. Models are inhabited by intensions, a
function $I$ sends terms to their intensions and functions $E_\alpha$ send
intensions to the extensions they determine. If an additional requirement
should be made that the $E_\alpha$ be injective, one essentially obtains
Henkin's general models, if, moreover, the
$E_{\langle\alpha_1\ldots\alpha_n\rangle}$ should be required to
be onto ${\mathcal P}(D_{\alpha_1}\times\cdots\times D_{\alpha_n})$,
standard models are obtained.

The logic contrasts with other approaches to (hyper-)intensionality in two
ways. Firstly, unlike other approaches, the aim is not to set up a new
logic, but to provide existing classical type theory with a wider class of
models in order to invalidate the axiom of Extensionality, which is unwanted
in many applications. Secondly, the logic is agnostic about
what intensions \emph{are}. To the latter question various answers have been
given but here we have only provided an abstract characterization of the
notion of intensionality. We have, in other words, focused on the logic
rather than on the ontology of intensions.

While the logic is a generalisation of classical type theory, not an
extension with new concepts, it turns out that there is a natural connection
with the usual notion of modality. Intensional models may have domains of
the propositional type $\langle\rangle$ that are not isomorphic with
$\{0,1\}$ and certain properties of objects in these domains can be
identified with possible worlds. Accessibility relations of various kinds
between such worlds are easily definable and modal box and diamond operators
can be obtained accordingly.

\section*{Acknowledgements} 
I wish to thank Nissim Francez for providing me with detailed comments and
for urging me to include linguistic applications. Prof.\ Roger Hindley very
kindly helped me with obtaining a copy of Takahashi's paper on
cut-elimination in type theory. The anonymous referee's highly welcome
comments led to several improvements, among which is the incorporation of
possible worlds semantics in section \ref{worlds}.

\bibliographystyle{plain}
\small
\newcommand{\SortNoop}[1]{}

\end{document}